\documentclass{amsart}
\date{}
\usepackage{amsmath}
\usepackage{amscd}
\usepackage{amssymb}
\usepackage{amsfonts}
\usepackage{latexsym}
\usepackage{verbatim}
\usepackage{epsfig}

\numberwithin{equation}{section}
\newtheorem{lemma}[equation]{Lemma}
\newtheorem{proposition}[equation]{Proposition}
\newtheorem{theorem}[equation]{Theorem}
\newtheorem{corollary}[equation]{Corollary}

\theoremstyle{remark}

\newtheorem{remark}[equation]{Remark}
\newtheorem{definition}[equation]{Definition}
\newtheorem{definitions}[equation]{Definitions}
\newtheorem{example}[equation]{Example}
\newtheorem{condition}[equation]{\textnormal{Condition}}
\newtheorem{questions}[equation]{Questions}

\newcommand{\eps}{\epsilon}
\newcommand{\fix}{\operatorname{Fix}}

\newcommand{\QED}{$\Box$}

\newenvironment{demo}
      {\medbreak\noindent{\sc Proof:}}
      {\hfill\QED\medbreak}

\begin{document}

\title
[Almost isomorphism ] {Almost isomorphism for countable state
Markov shifts}

\dedicatory{Dedicated to Peter Walters and Klaus Schmidt, on the
occasion of their sixtieth birthdays.}
\begin{abstract}
Countable state Markov shifts are a natural generalization of the
well-known subshifts of finite type. They are the subject of
current research both for their own sake and as models for smooth
dynamical systems. In this paper, we investigate their almost
isomorphism and entropy conjugacy and obtain a complete
classification for the especially important class of strongly
positive recurrent Markov shifts. This gives a complete
classification up to entropy-conjugacy of  the natural
extensions of smooth entropy-expanding maps, e.g., $C^\infty$
smooth interval maps with non-zero topological entropy.
\end{abstract}

\author{Mike Boyle}
\address{Mike Boyle\\Department of Mathematics\\
        University of Maryland\\
        College Park, MD 20742-4015\\
        U.S.A.}
\email{mmb@math.umd.edu} \urladdr{www.math.umd.edu/$\sim$mmb}
\author{Jerome Buzzi}
\address{Jerome Buzzi\\Centre de Math\'ematiques\\
Ecole polytechnique\\91128 Palaiseau Cedex\\France}
\email{buzzi@math.polytechnique.fr}
\urladdr{www.geocities.com/jrbuzzi/}
\author{Ricardo Gomez}
\address{Ricardo Gomez\\ Instituto de Matematicas Area de la Investigacion Cientifica\\
Circuito Exterior\\ Ciudad Universitaria\\ DF 04510\\ Mexico }
\email{rgomez@math.unam.mx}
\urladdr{www.math.unam.mx/$\sim$rgomez}

\thanks{The research of J. Buzzi was also supported by the Center for
Dynamical Systems of Penn State University and the kind
hospitality of the University of Maryland.}

\keywords{entropy; countable state Markov shift; topological
Markov chain; almost isomorphism; entropy conjugacy; magic word;
strong positive recurrence;  Artin-Mazur zeta function;
 smooth ergodic theory.}

\subjclass[2000]{Primary: 37B10; Secondary: 37B40, 37C99, 37D35}


\maketitle

\tableofcontents

\section{Introduction}

In this paper, ``Markov shift'' means a countable (possibly
finite) state irreducible Markov shift. We use the same symbol to
denote the domain of a Markov shift and the shift map on this
domain. There are several characterizations (\ref{sprprop}) of the
class of {\it strongly positive recurrent} (SPR)  Markov shifts;
this is  the  class of Markov shifts which  most resemble finite
state Markov shifts.

A map $\varphi  :S\to T$ between Markov shifts is a {\it one-block
code} if there is a function $\Phi $ from the symbol set of $S$
into the symbol set of $T$ such that
 $(\varphi  x)_n = \Phi  (x_n)$, for all $x$ and $n$.
A $T$-word $W$ (of length $|W|$) is a {\it magic word} for such a
map $\varphi $ if the following hold.
\begin{enumerate}
\item \label{mw1} If
 $y\in T$ and
$\{n\in \mathbb Z: y[n,n+|W|-1]=W\}$ is unbounded above and
unbounded below, then $y$ has a preimage under $\varphi $. \item
\label{mw2} There is an integer  $I$ such that whenever $C$ is a
$T$-word and two points $x$ and $x'$ of $S$ satisfy $(\varphi
x)[0,2|W|+|C|-1]=WCW= (\varphi  x')[0,2|W|+|C|-1]$, then $x[I,
I+|W|+|C|-1]= x'[I, I+|W|+|C|-1]$.
\end{enumerate}
(In the constructions of this paper, the integer $I$ of the last
condition will be zero.)
It follows from (\ref{mw2}) that the
preimage in (\ref{mw1}) is unique.

We define two Markov shifts $S$ and $T$ to be {\it almost
isomorphic} if there exist another Markov shift $R$ and injective
one-block codes $ R\to S$, $R\to T$ each of which has a
  magic
word  (recall that all Markov shifts are understood to be
irreducible). An almost isomorphism will induce a shift-commuting,
Borel bimeasurable bijection between the images of the one-block
codes, and thereby will induce isomorphisms of measurable systems
$(S,\mu )\to (T,\nu )$ for shift-invariant Borel probabilities
$\mu , \nu$
 which assign these images measure 1.
 The collection of such measures will include the
ergodic measures which give positive measure to every nonempty
open set,  and in the SPR case will include all ergodic measures
with entropy sufficiently close to the topological entropy
(defined (\ref{entropydefinition}) as the supremum of the
measure-theoretic entropies with respect to invariant Borel
probabilities). These isomorphisms will be finitary
(homeomorphisms between measure-one sets), with exponentially fast
coding time for exponentially recurrent ergodic measures (such as
the unique measure of maximal entropy, in the SPR case).

Our main result (\ref{maintheorem}) is that SPR Markov shifts
 are almost isomorphic if and only if they
have the same entropy and period. This can be viewed as an
 analogue of the Adler-Marcus classification of
irreducible shifts of finite type up to almost topological
conjugacy by entropy and period (\ref{adlermarcus}). The
classification beyond SPR cannot possibly have the same simplicity
 (\ref{ulfremark}).
We also give a sufficient condition for almost isomorphism of not
necessarily SPR Markov shifts (\ref{sufficienttheorem}).

Two systems are {\it entropy conjugate} if they are Borel
conjugate after restriction to sets which have full measure for
all ergodic measures with entropy near the topological entropy
(see Definition \ref{entcondef}). Various smooth, piecewise smooth
and symbolic systems  are known  (Theorem \ref{classlist}) to have
natural extensions which are entropy conjugate to a finite union
of SPR Markov shifts. Consequently, our main result provides
simple
 invariants which classify the natural extensions of
these systems up to entropy conjugacy (Theorem
\ref{piecewiseclassif}).

 One motivation of our paper is recent interest in the
thermodynamic formalism for countable state Markov shifts
\cite{BS,FFY,G1,GP,GS,MU1,Sarig0,Sarig,Sarig2,Sarig3}. In
 a companion paper \cite{BBG},
 we describe a reasonable class of Borel potentials which behaves
well  under almost isomorphism. We are also motivated by the use
of countable state Markov shifts to code some partially or
piecewise hyperbolic systems
(\cite{SIM,EE,QFT,H1,H2,H3,K1,K2,Tak,Tak2,MU2,Young} and
successors of \cite{Young}); investigations of coding relations
among Markov shifts (see \cite{DF,DF2,DF3,DF4,FF} and their
references); and a longstanding finite state coding problem for
Markov measures not of maximal entropy (\ref{othermarkov}).

This paper is dedicated to
Klaus Schmidt and Peter Walters, on the occasion of their sixtieth
birthdays.
Words  from the poet E.E.Cummings
\cite{EEC} suit them:
\begin{quote} {\sl
septembering arms
of year extend \\
less humbly wealth to fool and friend }
\end{quote}

\section{Definitions and background}

 Let $A$ be a square matrix with nonnegative integer
entries, with rows and columns indexed by a finite or countable
set $\mathcal I$. View $A$ as the adjacency matrix of a directed
graph with some set of edges $\mathcal E_A$. Give $\mathcal E_A$
the discrete topology; give $(\mathcal E_A)^{\mathbb Z}$ the
product topology; let $S_A$ denote the set of doubly infinite
sequences on $\mathcal E_A$ corresponding to walks through the
graph; and let $S_A$ have the topology it inherits as a subspace
of
 $(\mathcal E_A)^{\mathbb Z}$. The shift map homeomorphism on $S_A$
is defined by sending a bisequence $(x_n)$ to the bisequence
$(y_n)$ such that $y_n=x_{n+1}$ for all $n$. Given a
point/bisequence $x$, we let $x[i,j]$ denote the word
$x_ix_{i+1}\dots x_j$. If $W$ is a path of $n$ edges (or a finite
sequence of $n$ symbols), then $|W|$ denotes $n$. We will use the
same symbol (e.g. $S_A$) for the shift map and its domain. $S_A$
is uniquely determined by $A$, up to the naming of the edges in
the graph;  in this paper, we can safely ignore this ambiguity,
and refer to ``the'' Markov shift defined by $A$.
 We use the edge
shift presentation rather than the vertex shift presentation, but
this paper could have been written entirely with the vertex shift
presentation.

The matrix $A$ is {\it irreducible} if for every row $i$ and
column $j$ there exists $n>0$ such that $A^n(i,j)>0$. (The term
``indecomposable'' as used in \cite{GS} is synonymous with our
``irreducible''.) Convention: in this paper, by a  {\it Markov
shift} we will mean a homeomorphism $S_A$ defined by a (finite or
countably infinite) {\it irreducible}
 matrix $A$ over $\mathbb Z_+$.

By a graph we will always mean a directed graph, with finitely or
countably infinitely many vertices  and edges. Given a vertex $v$
in a graph, a first return loop to $v$ is a path of edges which
begins and ends at $v$ and otherwise does not visit $v$. The
number of edges in a path $p$ is its {\it length},  $|p|$. A {\it
loop graph} (in \cite{GS}, a {\it petal graph}) is a graph
$\mathcal G$ with a distinguished vertex $v$, such that $\mathcal
G$ is the union of the first return loops to $v$, and every vertex
except $v$ lies on a unique first return loop. From here, for
brevity by a {\it loop} in a loop graph we will always mean a
first return loop to the distinguished vertex.

Let $f$ be a  power series  $\sum_{n=1}^{\infty }f_nz^n$ with
nonnegative integer coefficients, i.e.  $f\in z\mathbb Z_+[[z]]$.
We will sometimes use the notation rad($f$) to denote the radius
 of convergence of $f$. A {\it loop graph} $\mathcal G_f$
for $f$ is a loop graph which for each $n$ contains exactly
 $f_n$ loops of length $n$. The loop graph is uniquely
determined by $f$ up to the naming of the edges and vertices.
Given
 $\mathcal G_f$,
the {\it loop shift} $\sigma_f$ is the Markov shift $S_A$ such
that $A$ is the adjacency matrix of a loop graph for $f$.

Suppose $S_A$ is  a Markov shift,  with associated graph $\mathcal
G(A)$,  and $v$ is a vertex in $\mathcal G(A)$. Let $f_n$ denote
the number of first return loops of length $n$ to $v$,   and
assume $f_n<\infty$ for all $n$ (a property which is independent
of the choice of vertex). Let $\sigma_f$ be the corresponding loop
shift. The {\it period} of $S_A$ (or of $\sigma_f$) is the g.c.d
of $\{n: f_n >0\}$, and does not depend on the choice of vertex.
$S_A$ is mixing if and only if the period is 1.

\begin{definition}\label{entropydefinition}
By the entropy $h(S_A)$($=h(\sigma_f)$), we will mean the
 Gurevich entropy
\cite{G0}, which is the supremum of the measure theoretic
entropies over invariant Borel probabilities.
This entropy equals $ \log (\varlimsup |r_n|^{1/n})$, where $r_n$
is the number of loops (not necessarily of first return) at an
arbitrary vertex, and
 the limsup here is a limit
when $S_A$ is mixing.
\end{definition}
Let $\log (\lambda )$ denote the entropy of $S_A$ and suppose
$0<\log(\lambda )< \infty$. Four basic classes can be described in
terms of the sequence $(f_n)$.
\begin{enumerate}
\item $S_A$ is transient if  $\sum f_n/\lambda^n< 1$ \item $S_A$
is recurrent if  $\sum f_n/\lambda^n =1 $ \item $S_A$ is positive
recurrent if $\sum f_n/\lambda^n =1 $ and $\sum n
f_n/\lambda^n<\infty $ \item $S_A$ is strongly positive recurrent
(SPR) if $\varlimsup |f_n|^{1/n} < \lambda $.
\end{enumerate}
(We are using names of matrix classes to describe the Markov
shifts they define.) Here  (4)$\implies$(3)$\implies$(2). The
properties (1)-(4) don't depend on the particular choice of vertex
\cite{GS}. The classes (1)-(3) correspond to the classical
Vere-Jones classification \cite{Vere-Jones1,Vere-Jones2} of
$\mathbb N \times \mathbb N$ nonnegative matrices. A finite
entropy Markov shift  is positive recurrent if and only if it has
a measure of maximal entropy $\mu$ (i.e., $h(S)=h(S,\mu )$), in
which case there is only one measure of maximal entropy \cite{G0}.

Markov shifts defined over a countable alphabet may differ
significantly from those which have a finite alphabet (subshifts
of finite type), see e.g. \cite{DF,GS,Ki,Sarig2,Sarig3}. The
search for well-behaved Markov shifts leads to the class of SPR
shifts. The notion of SPR  already appeared in Vere-Jones' work
\cite{Vere-Jones2} in the special case of stochastic matrices as a
necessary and sufficient condition for some exponential
convergence. Following work of Salama (see
\cite{ruette,salama1,salama2}) on conditions for the strict
decrease of entropy in proper subsystems, the class of SPR shifts
was introduced independently by  Ulf Fiebig \cite{UF} and by
Gurevich \cite{G1}, who introduced the class of ``stably
recurrent'' matrices, which contain those defining what we  call
the SPR Markov shifts. In \cite{GS}, this class is developed
further as the fundamental class of ``stable positive'' matrices
\cite{GS}.  We view ``SPR'' as also abbreviating ``stable positive
recurrent''.

We will record some of the conditions on a Markov shift which are
equivalent to SPR. For one, we first recall the definition of
exponential recurrence while remarking some equivalences.

\begin{remark} \label{exprec}
Suppose $(S,\mu )$ is an ergodic automorphism of a
topological space with a Borel probability measure. Let $V$ be a
measurable set such that $\mu V>0$. Let $r_V$ be the return-time
function on $V$, $r_V(x)=\min\{n>0:S^nx\in V \}$. By  ergodicity
the system $(S,\mu )$ is isomorphic to a tower over $V$ with
return time function $r_V$. It is not difficult to check  that the
following are equivalent:
\begin{enumerate}
\item $(S,\mu )$ is exponentially recurrent: i.e. for every open
set $V$ with $\mu V>0$,
\[
\varlimsup_n \big(\mu \{x\in V : r_V(x) \geq n\}\big)^{1/n} < 1 \
.
\]
\item $(S,\mu )$ is exponentially filling: i.e. for every open set
$V$ with $\mu V>0$,
\[
\varlimsup_n \big(\mu \{x \in S :x \notin \cup_{k=1}^n S^{-k}V
\}\big)^{1/n} < 1\ .
\]
\item $(S^{-1},\mu )$ is exponentially recurrent.
\end{enumerate}
\end{remark}

\begin{proposition}\label{sprprop}
The following are equivalent conditions on a Markov shift $S$.
 \begin{enumerate}
\item \label{-spr} $S$ is SPR. \item Removing any edge from the
graph defining $S$ strictly lowers the  entropy (\cite{UF},
\cite{salama2},\cite[Remark 3.16]{GS}; see  \cite{ruette}).
  \item
Some (equivalently every)
 local zeta function (counting the number of fixed points $x\in
[v]$ for some vertex $v$)   has a non-trivial meromorphic
extension (see \cite{GS} and Section \ref{sectionzeta}).

 \item \label{-exprec}
$S$ has a measure of maximal entropy $\mu$, and $(S,\mu )$ is
exponentially recurrent.
 \end{enumerate}
\end{proposition}

\begin{proof}
We will prove (\ref{-exprec})$\iff$(\ref{-spr}). Note
that both conditions imply that there exists a measure of maximum
entropy and recall that such a measure gives positive measure to
every non-empty open subset \cite{G0}. We may and do assume that
we have such a measure $\mu$.

Let  $V_I$ be the open set of points $x$ such that the edge $x_0$
begins at a given vertex $I$. Let $f_n$ be the number of first
return loops to $I$. By the specific form for $\mu$ \cite{Ki,G0},
there is a constant $c_I>0$ such that, for any first return loop
$W_0\dots W_{n-1}$, the set
 $[W]:=\{x: x_0\cdots x_{n -1}=W\}$ has $\mu$-measure
$c_I(1/\lambda )^{n}$. Thus $\mu (S\setminus \cup_{n=1}^NS^nV_I)$
goes to zero exponentially fast, i.e. $V_I$ is exponentially
filling in $(S,\mu)$, if and only if $\sum_{n=N+1}^{\infty}
c_I(n-N)(1/\lambda )^{n}$ goes to zero exponentially fast.  But
the latter is equivalent to the SPR property.

 To conclude, it is enough to see that SPR implies not
only exponential filling for the sets $V_I$ for any vertex $I$,
but for all non-empty open sets $V$. Any such $V$ contains a
cylinder set $C= \{x:x[-k,k]=W\}$ for some word $W$, and $W$
becomes a vertex in a higher block presentation $S'$ of $S$.
Because $S'$ is again SPR \cite[Propostion 2.12]{GS}, $C$ is
exponentially filling in $S'$, and therefore so is $V$ in $S$.
\end{proof}

For more background, see
 \cite{GS,Ki,ruette} and their
references for Markov shifts and nonnegative matrices;
 \cite{Walters-book}
for entropy and ergodic theory; and \cite{Ki,LM} for symbolic
dynamics.

\section{Magic words, almost isomorphism and entropy-conjugacy}

We assume the definitions already given in the Introduction. Given
a subset $K$ of a Markov shift $S$, we let $\mathcal M(K)$ denote
the set of shift-invariant Borel probabilities $\mu$ on $S$ such
that $\mu K=1$. Let $K'$ be a subset of another Markov shift $S'$.
Given a shift-commuting Borel automorphism $\gamma :K\to K'$, we
may use the same symbol $\gamma$ for related maps such as $\gamma
: \mathcal M(K)\to \mathcal M(K')$. Recall a measure on a space
has {\it full support} if it is nonzero on every nonempty open
set. If $\mu$ is a shift invariant measure on $S$, then it defines
a measure-preserving system $(S,\mu )$. An isomorphism of measure
preserving shift systems $(S,\mu )\to (S',\mu ')$ is a
bimeasurable, shift-commuting,
 measure-preserving bijection from a subset $K$ of $S$
to a subset $K'$ of $S'$, where $\mu K=\mu 'K'=1$. Such an
isomorphism is {\it finitary} if
 $K$ and $K'$ can be chosen so that
the map $K\to K'$ is a homeomorphism, with respect to the
topologies $K$ and $K'$ inherit as subsets of the topological
spaces $S$ and $S'$. In this case, for $\mu$-almost all $x$, there
exists  a minimal nonnegative integer $n=n(x)$ such that for
$\mu$-almost all points $x'$ with $x'[-n,n]=x[-n,n]$, we have
$(\varphi  x)_0=(\varphi  x')_0$. (In other words, $\mu$-almost
surely $\varphi  $ is a variable length block code.) We  then
define the {\it expected coding time} of $\varphi $  to be $\int_S
n(x) d\mu $. We say $\varphi $ has {\it exponentially fast} coding
if $\varlimsup_k |\mu \{x:n(x)\geq k\}|^{1/k}<1$; this property
implies the expected coding time is finite.
 We say systems $(S,\mu )$ and $(S',\mu ')$ are
{\it finitarily isomorphic} if there exists a finitary isomorphism
between them. We say following \cite{S2,S3} that a measurable
isomorphism is {\it hyperbolic structure preserving} if on sets of
full measure it respects the stable and unstable relations (here,
the relations of being forwardly/backwardly asymptotic under the
shift).

For $S,S'$ Markov shifts with given invariant measures $\mu\mapsto
\mu'$, a {\it magic word isomorphism} is a  finitary isomorphism
$\gamma :(S,\mu)\to (S',\mu ')$ such that both $\gamma $ and
$\gamma^{-1}$ have a magic word. The definition of magic word here
is the same as in the introduction, with the modifications that
conditions (1) and (2) are required to hold only on sets of full
measure, and the set of points seeing the magic word $W$ in
infinitely many positive and negative coordinates is required to
have full measure. The latter condition is automatic for an
ergodic measure assigning positive measure to the set $\{ x:
x[0,|W|-1]=W\}$.
(We remark that the definition of magic word isomorphism in
\cite{G} is slightly different; a magic word here for $\gamma$
is a magic word in \cite{G} for $\gamma^{-1}$.)

\begin{proposition} \label{expspeed}
Suppose $\gamma : (S, \mu )\to (S',\mu ')$ is a magic word
isomorphism between  Markov shifts endowed with invariant and
ergodic probability measures. Then $\gamma$ is a hyperbolic
structure preserving isomorphism. If the systems are exponentially
recurrent, then $\gamma $ codes exponentially fast.
\end{proposition}

\begin{proof}
With probability one, points will see a magic word infinitely
often in positive and in negative coordinates. Therefore, the left
tail $x(-\infty ,0]$ of a bisequence $x$ will code some left tail
$y(-\infty ,m]$ of its image, and likewise for right tails. This
proves the a.e. preservation of hyperbolic structure.

Now suppose the systems are exponentially recurrent. Let $V$ be an
$S$-word which is a
 magic word for $\gamma^{-1}$.
Without loss of generality, let the index  $I$ used in the
definition of a magic word in the Introduction be zero. Let  $J$
be the length of $V$. For $n>J$, let $E_n$ be the set of
$S$-points $x$ such that $V$ is not a subword  of $x[-n,-1]$, and
let $F_n$ be the set of $S$-points $x$ such that $V$ is not a
subword of $x[1,n]$. Then $x[-n,n]$ codes $(\varphi  x)_0$ except
for a subset of $E_n \cup F_n$. Let $V$ also denote the set $\{x:
x[0,J-1]=V\}$. We have
\[
F_n = \{x: x\notin \cup_{k=1}^{n-J} S^{-k}(V)\} \quad , \quad E_n
= \{x: x\notin  \cup_{k=1}^{n-J} S^k(S^{J}V)\}
\]
and it follows from Remark \ref{exprec} that
\[
\varlimsup_n \big(\mu (E_n \cup F_n)\big)^{1/n} \leq \varlimsup_n
\big(\mu E_n + \mu F_n\big)^{1/n} <1 \ .
\]
\end{proof}

\begin{definitions}\label{measdef}
By a {\it measurable  system} we will mean a measurable map
$T:X\to X$. In this paper, the underlying $\sigma$-algebra will be
the Borel $\sigma$-algebra naturally associated to $X$, and the
entropy $h(T)$ is
 defined as the supremum of the measure theoretic
 entropies with respect to
$T$-invariant Borel probabilities.
For a measurable system $T:X\to X$,   a subset $N\subset X$ is
called {\it entropy-negligible} \cite{SIM} if there is $h< h(T)$
such that $\mu(N)=0$ for every ergodic invariant
probability measure $\mu$ with $h(T,\mu)> h$.
\end{definitions}

\begin{definition}\label{entcondef}
Two measurable systems $T:X\to X$ and $S:Y\to Y$ are {\it
entropy-conjugate}
 if there exist entropy-negligible subsets
$X_0\subset X$ and $Y_0\subset Y$ and a bimeasurable bijection
$\gamma : X\setminus X_0 \to Y\setminus Y_0$ such that $S\gamma =
\gamma T$ for all $x\in X\setminus X_0 $. Such a map $\gamma$ is
called an {\it entropy-conjugacy}.
\end{definition}

\begin{proposition}\label{magicprop}
Suppose $S$ and $T$ are almost isomorphic Markov shifts. Then
$h(S)=h(T)$, and there are Borel subsets $K$ and $K'$ of $S$ and
$T$, and a shift-commuting Borel-measurable bijection $\gamma
:K\to K'$, such that the following hold.
\begin{enumerate}
\item \label{ai1} $K$ and $K'$ are residual subsets of $S$ and $T$
(contain dense $G_{\delta}$ sets). \item \label{ai2} The map
$\gamma$ induces a bijection
 $\mathcal M(K)\to \mathcal M(K')$ ($\mu \mapsto \mu '$, say)
such that for each such pair $\mu,\mu'$
 the map $\gamma$ induces an
 isomorphism $\gamma :(S,\mu )\to (T, \mu ')$,
which is a magic word isomorphism when $\mu$ and $\mu'$
have full support.
 \item
\label{ai3} $\mathcal M(K)$ and $\mathcal M(K')$
 contain all ergodic shift-invariant Borel
probabilities on $S$ and $T$ with full support, and these
correspond under
$\gamma$. \item \label{ai5} If $S$ is SPR, then so is $T$, and
$\gamma $ is an  entropy-conjugacy from $S$ to $T$.
\end{enumerate}
\end{proposition}

\begin{demo}
Suppose $S$ and $T$ are almost isomorphic Markov shifts, i.e.,
suppose we have another Markov shift $R$ and injective one-block
codes $\varphi :R\to S$ and  $\psi: R\to T$, with magic
words $W$ and $W'$ respectively. We assume for simplicity that
$I=0$ in their definition.

Let $W_*$, resp. $W_*'$, be some $R$-word projecting to $W$, resp.
$W'$. Pick some other $R$-words $C_*,D_*,E_*$ such that
$w_*=W_*C_*W'_*D_*W'_*E_*W_*$ is also an $R$-word. Let $w$ be the
$S$-word below this  $R$-word. Each $x\in
\operatorname{Fix}(S^n)$ such that $x[0,|w|-1]=w$ has a unique
preimage $\varphi^{-1}(x)\in\operatorname{Fix}(R^n)$ and it
satisfies $\varphi^{-1}(x)[0,|w|-|W_*|-1]=W_*C_*W'_*D_*W'_*E_*$.
Therefore
$\psi\varphi^{-1}(x)[|W_*C_*|,|W_*C_*W'_*D_*W'_*|-1]=W'CW'$ for
some $T$-word $C$.
$\psi\varphi^{-1}(x)\in\operatorname{Fix}(T^n)$. Moreover the map
$\psi$ is injective on $\operatorname{Fix}(R^n)\cap[w_*]$ as $w_*$
projects in $T$ to a magic word.

Thus the number of $n$-periodic sequences in $S$ which starts with
$w$ is a lower bound for the number of $n$-periodic sequences in
$T$ which starts with $W'CW'$. The growth rates of these numbers
are the entropies.  Thus $h(S)\leq h(T)$. By symmetry,
$h(S)=h(T)$.

(\ref{ai1}) Let $K$ be the image of $\varphi $ in $S$ and let $K'$
be the image of $\psi$ in $T$. Let $\gamma :K\to K'$ be the
bijection $\gamma = \psi \varphi ^{-1}$. For any $n$, the set
$K_n(W)$ of points in $S$ which see $W$ in at least $n$ negative
coordinates and $n$ positive coordinates is a dense open set, and
$K$ contains the intersection $K(W)$ of the $K_n(W)$. Thus $K$ is
a residual subset in $S$, as is $K'$ in $T$.

(\ref{ai2}) With respect to the measure 1 sets  $K(W)\cap
\gamma^{-1}K(W')$ and $\gamma K(W) \cap K(W')$, the map $\gamma$
is a magic word isomorphism  between corresponding
measures with full support.

(\ref{ai3}) Let $C_W=\{x: x[0,|W|-1]=W\}$. Let $\mu$ be ergodic
with full support; then $\mu C_W>0$, and by ergodicity
$\mu (K(W))=1$, so $\mu(K)=1$ and $\mu\in\mathcal M(K)$. Suppose
$U$ is an open set in $S$;  then $U\cap K$ is open in $K$ and $\mu
U=\mu(K\cap U)=\mu '(\gamma (K\cap U))$. Modulo null sets,
non-empty open sets of $K'$ contain such non-empty sets
$\gamma(K\cap U)$. It follows that $\mu'=\gamma\mu$ is a fully
supported measure when $\mu$ is. The converse is shown similarly.

(\ref{ai5}) $\gamma$  takes  the measure of maximal entropy $\mu$
of $S$ (which has full support \cite{G0}) to a measure $\mu'$ of
$T$ which is isomorphic and therefore of maximal entropy (recall
$h(S)=h(T)$). Because $S$ is SPR, by Proposition \ref{sprprop}
$(S,\mu )$ is exponentially recurrent. By \cite{RudMarkov},
exponential recurrence is invariant under finitary isomorphism, so
$\mu '$ is an exponentially recurrent measure of maximal entropy
for $T$. We conclude from
  Proposition   \ref{sprprop} that
$T$ is SPR.

 By Proposition \ref{sprprop}, the entropies of measures
such that $\mu(C_W)=0$ are bounded away from $h(S)$, hence
$\mathcal M(K)$ contains all ergodic and invariant measures with
entropy close to $h(S)$. The same holds for $\mathcal M(K')$. By
(\ref{ai2}), $\gamma$ is an entropy-conjugacy.
\end{demo}

\begin{remark}\label{hc-ai-remark}
Suppose $S,S'$ are SPR Markov shifts; $\mu$,$\mu'$ are ergodic
measures with full support; and $\gamma : (S,\mu )\to (S,\mu ')$
is a magic word isomorphism. The previous proof shows  $\gamma$
must be the restriction of
 an entropy-conjugacy between $S$ and $S'$.
\end{remark}

\begin{remark} \label{ulfremark}
For positive recurrent shifts which are not SPR, the work of Ulf
Fiebig  \cite{UF2} shows that recurrent rate invariants of
finitary isomorphism are much finer and richer. These also give
invariants of almost isomorphism.
\end{remark}

\begin{proposition} \label{trueprop} The following are true.
\begin{enumerate}
 \item Suppose $R\to S$ and $S\to T$ are injective one-block codes
between Markov shifts, each of which has a magic word. Then their
composition is an injective one-block code with a magic word.
 \item Suppose $\sigma_f$ is the loop system built from
first-return loops to a vertex for a Markov shift $S_A$. Then the
natural inclusion $\sigma_f\to S_A$ is a one-block code with a
magic word.
 \item If $R\to S$, $R\to T$ is an almost
isomorphism between $S$ and $T$, and $\sigma_f$ is a
first-return-to-a-vertex loop system for $R$, then there is also
an almost isomorphism $\sigma_f\to S$, $\sigma_f\to T$.
 \item Almost isomorphism defines an equivalence
relation on Markov shifts.
\end{enumerate}
\end{proposition}

\begin{proof}
We leave the proof of Proposition \ref{trueprop} as an exercise.
 For the last claim (which we do not use in
our proofs),  the nontrivial step is to show, given Markov shifts
$R_1,R_2,T$ and injective one block codes with magic words
$\varphi  :R_1\to T$, $\varphi _2 : R_2\to T$, that there is a
Markov shift $S$ and injective one block codes $p_1: S\to R_1$,
$p_2:S\to R_2$. For this one can use a standard fiber product
construction
 \cite{Ki,LM}, as follows.
Let $S'$ be the subset of $R_1\times R_2$ which consists of points
$(x,y)$ such that $\varphi _1x=\varphi _2y$, and define $p_1:
(x,y)\mapsto x$ and $p_2: (x,y)\mapsto y$. Let $S$ be the unique
maximal irreducible Markov shift contained in $S'$ which has dense
image in $R_1$ and in $R_2$. It can be verified that $S,p_1,p_2$
meet the requirements.
\end{proof}

Within a class of almost isomorphic Markov shifts, one can freely
localize to loop shifts and still capture all ergodic phenomena
with full support. However, if a Markov shift $S$ is not SPR, the
localization may miss many periodic points and measures of large
entropy, as the next example shows.

\begin{example}
Suppose $\sigma_f$ is not SPR. Then \cite{ruette} there is a loop
$\ell$ of $\sigma_f$ such that the subsystem $T$ which misses
$\ell$ has $h(T)=h(S)=\log \lambda$. Pick a vertex on $\ell$
(other than the base vertex for $\sigma_f$) and let $\sigma_g$ be
the loop shift based at this vertex. We have the obvious injective
block code $\sigma_g\to \sigma_f$ with a magic word. However,
 the supremum of
the entropies of ergodic measures
 assigning measure one to the complement of the
image is $\log\lambda$. Therefore, the inclusion
 $\sigma_g\to \sigma_f$ is not an entropy-conjugacy.
\end{example}

\section{Zeta functions}\label{sectionzeta}

Suppose $T$ is a bijection of a set, $T:X\to X$, such that for all
$n>0 $ the cardinality of the fixed point set of $T^n$ (i.e.,
$|\textnormal{Fix}(T^n)|$) is finite. Then the Artin-Mazur zeta
function of $T$ is defined as
\[
    \zeta_T(z) := \exp \sum_{n\geq1}
|\textnormal{Fix}(T^n)|
 \frac{z^n}{n} \  .
\]
Let $\mathcal O_n(T)$ denote the set of orbits of cardinality $n$
of $T$. We have the well known product formula
\begin{equation}
   \zeta_T(z)
= \Big( {\prod_{n=1}^{\infty}(1-z^n)^{|\mathcal
O_n(T)|}}\Big)^{-1}
 \ .
\end{equation}
This formula is easily verified when $T$ contains just one finite
orbit, and then holds at the level of formal power series by
taking countable products. (Note, the zeta function and product
formula make perfectly good sense at the level of formal power
series even when $\varlimsup |\textnormal{Fix}(T^n)|^{1/n}=\infty
$.) For a loop shift $\sigma_f$, with $f(z)=\sum_{n=1}^{\infty
}f_nz^n$, we have
\begin{equation}\label{loopzeta}
\zeta_{\sigma_f}(z) = \frac 1{1-f(z)}
 \ .
\end{equation}
We do not know the earliest reference for (\ref{loopzeta}); it
appears in \cite{Tak2} and in the case $f(z)$ is a polynomial it
is a special case of a result proved in
 the ``bigamy'' paper  \cite{BGMY}.
For a detailed proof of a generalization of
 (\ref{loopzeta}), see
\cite[Theorem 8.2]{GS}. We remark, (\ref{loopzeta}) in general
follows from the case that $f(z)$ is a polynomial,
 because coefficients
of degree at most $n$ in the series on either side of
(\ref{loopzeta}) match what would be computed for the loop shift
for the truncated series $f_1z+f_2z^2+\cdots +f_nz^n$.

A Markov shift $T$ may fail to have a well defined zeta function
in that $|\textnormal{Fix}(T^n)|$ may be infinite for some or all
$n$.  To avoid this problem at least when the entropy is finite,
we use, as in \cite{GS,Sarig}, a local zeta function, i.e., the
zeta function of the loop shift given by the return loops to a
given vertex $s$. For a finite entropy Markov shift, for every $n$
there are only finitely many first return loops to $s$ of length
$n$, so the local zeta function must be well defined.

\begin{proposition}\label{gsprop}
Let $\sigma_f$ be a  mixing SPR loop shift with finite positive
entropy $\log (\lambda )$. Let $\zeta (z)$ denote its zeta
function $\zeta_{\sigma_f}(z)$. Then $\zeta (z)$ is holomorphic in
$|z|< 1/\lambda$ and has a meromorphic extension to a larger disk
 (all our disks are centered at zero). More precisely,
$(1-\lambda z)\zeta(z)$ is holomorphic and nonzero in some disk
$|z|<r$ with $r>1/\lambda$.
\end{proposition}

Proposition \ref{gsprop} is  in the treatise of
 Gurevich and Savchenko \cite[Proposition 9.2]{GS}. We give
a proof for completeness.

\begin{demo}
By definition of SPR, $f$ is holomorphic on a disk of radius
strictly greater than $\lambda^{-1}$. Its power series has no
constant term, hence it is not constantly equal to $1$. It then
follows from (\ref{loopzeta}) that $\zeta (z)$ has a meromorphic
extension to this disk with poles exactly at points $z$ where
$f(z)=1$. Let $\zeta (z)$ also denote this meromorphic extension.

Now $f(\lambda^{-1})=1$ as  $\sigma_f$ is not transient.
Hence $\lambda^{-1}$ is a pole of $\zeta (z)$. As $f(z)$ is a
power series with non-negative coefficients not all zero,
$f'(\lambda^{-1})>0$. In particular,  $f'(\lambda^{-1})$ is
non-zero and the pole at $\lambda^{-1}$ is, as claimed, a simple
pole.

The mixing assumption implies that there are no other poles on the
circle $|z|=\lambda^{-1}$. Indeed, let $|z|=\lambda^{-1}$. We have
\[
    \left|\sum_{n\geq1} f_nz^n\right| \leq \sum_{n\geq1} f_n|z^n| = 1
\]
and the inequality is an equality iff $z^n=|z^n|=\lambda^{-n}$ for
all $n\geq1$ such that $f_n\ne0$. Such a $z$ must be $\omega
\lambda^{-1}$ with $\omega$ a root of unity of some order $q\geq1$
such that
 $f_n\ne 0\implies q|n$. Therefore $q$ divides the period of $\sigma_f$. But
$\sigma_f$ is mixing and therefore this period is $1$. In
particular, $q=1$ and $z=\lambda^{-1}$.  By compactness,
$z=\lambda^{-1}$ is the only pole of $\zeta (z)$ on some disk of
radius $>\lambda^{-1}$.
\end{demo}

The following corollary will be essential to the proof of our main
result (\ref{maintheorem}).

\begin{corollary} \label{coro-On}
Let $\sigma_f$ and $\sigma_g$ be mixing SPR loop shifts of equal
finite entropy $\log (\lambda )>0$. Then
\[
\varlimsup_n \Big| |\mathcal O_n(\sigma_f)| - |\mathcal
O_n(\sigma_g)| \Big|^{1/n}
 < \lambda \ .
\]
\end{corollary}

\begin{demo}
Set $a_n= |\fix ((\sigma_f)^n)| - |\fix ((\sigma_g)^n)| $ and
 $b_n= |\mathcal O_n(\sigma_f)| -
 |\mathcal O_n(\sigma_g)|$. Without loss of generality,
we assume $b_n\neq 0$ for infinitely many $n$, so $\varlimsup
|b_n|^{1/n}\geq 1$ and $\varlimsup |a_n|^{1/n}\geq 1$.
  Because $a_n = \sum_{k|n} kb_k$, and by M\"{o}bius inversion
 $n b_n = \sum_{k|n} \mu(n/k)a_k$,
 it then
follows that $\varlimsup |b_n|^{1/n} = \varlimsup |a_n|^{1/n}$.

By Proposition \ref{gsprop},
\[
\frac{\zeta_{\sigma_f}(z)}{\zeta_{\sigma_g}(z)} = \exp
\sum_{n\geq1} \Big( |\fix((\sigma_f)^n)| - |\fix((\sigma_g)^n)|
\Big) \frac{z^n}{n}
\]
 is a
holomorphic function on a disk of radius $r>1/\lambda$, and
consequently so is $\sum (a_n/n)z^n$. Therefore
$\varlimsup |a_n|^{1/n} < \lambda $, which implies
 $\varlimsup |b_n|^{1/n} < \lambda $.

\end{demo}

\section{The Loops Lemma}

Given a power series $k$ with zero constant term, we use the
standard notation $k^*$ to denote the power series $ 1 + k +
k^2+\dots =1/(1-k)$. Similarly, given a set $K$ of words,
containing only finitely many words of any fixed length, we let
$K^*$ denote the set of words which are concatenations of words in
$K$, and we include in $K^*$ the empty word, which we regard as an
identity element for concatenation.

The construction of the next lemma is lifted from \cite{G}.
\begin{lemma}\label{lemmazero}
Suppose $f=h+k$ where $h$ and $k$ are formal power series in
$z\mathbb Z_+[[z]]$ which are not identically zero.
 Then $1-hk^*= (1-f)/(1-k)$, and
there is an injective one-block code $\psi$ from $\sigma_{hk^*}$
to $\sigma_{f}$ which has a
 magic word. If additionally $h(\sigma_k)<h(\sigma_f)$,
then $\psi$ is also an entropy-conjugacy.
\end{lemma}
\begin{demo}
 The equality is an
application of the geometric series.

We have $h=\sum_{n=1}^{\infty} h_nz^n$ and
 $k=\sum_{n=1}^{\infty} k_nz^n$.
Let $\mathcal G_f$ be a loop graph for $f=h+k$, with $\mathcal L$
the set of (first return) loops of $\mathcal G$. Choose subsets
$\mathcal H$ and $\mathcal K$ of   $\mathcal L$ such that
$\mathcal L$ is the disjoint union of $\mathcal H$ and $\mathcal
K$; and for every $n$, the number of loops of length $n$  in
$\mathcal H$ is $h_n$; and for every $n$, the number of loops of
length $n$ in $\mathcal K$ is $k_n$.
 Let  $\mathcal L'$
be the loop set of a loop graph $\mathcal G_{hk^*}$ for the series
$hk^* = h(1+k+k^2+\cdots )$.
 Choose a bijection $\beta : \mathcal L'\to \mathcal H\mathcal K^*$
respecting word length and let $\beta $ define a labeling of
$\mathcal L'$, and thereby a one-block code $\psi$ from
$\sigma_{hk^*}$ to $\sigma_{f}$. Each loop $l$ of $\mathcal L'$
has a label $\widehat l$ which has a nonempty prefix which is an
element of
 $\mathcal H$, and no word of $\mathcal H$ occurs in
$\widehat l$ except as a prefix. Consequently any word of
$\mathcal H$ is a magic word for $\psi$.

 It follows in particular that $\psi$ is an isomorphism
w.r.t. any ergodic and invariant measure, unless the measure
avoids all words of $\mathcal H$ and thus lives on $\sigma_k$. The
entropy of such a measure is bounded by $h(\sigma_k)$, so
$\varphi$ is an entropy-conjugacy if $h(\sigma_k)<h(\sigma_f)$.
\end{demo}

The rest of this
 section is devoted to the following lemma, which records the
key features of our main construction.

\begin{lemma}[Loops Lemma] \label{loopslemma}
Let $\sigma_f$ be a mixing loop shift,
$f=\sum_{n=1}^{\infty}f_nz^n$, with
$0<h(\sigma_f)=\textnormal{log}(\lambda )\leq\infty$. Let
$(r_k)_{k=1}^{\infty}$ be a nondecreasing sequence of positive
integers\footnote{$r_k$ will be the length of the orbit to be
deleted at stage $k$.}. Let $R_n=\# \{k:r_k=n\}$. Set $f^{<1>}=f$.
Given
 $f^{<k>}$, define the series
$f^{<k+1>}=\sum_{n=1}^{\infty}f_n^{<k+1>}z^n$ by the equation
\[
f^{<k+1>}(z)= (f^{<k>}(z)-z^{r_k})\sum_{n=0}^{\infty}z^{nr_k}
\]
or equivalently
\[
1-f^{<k+1>}(z)= \frac{1-f^{<k>}(z)}{(1 -z^{r_k})}\ .
\]
Assume the pair $(f,(r_k)_{k=1}^{\infty})$ satisfies the
positivity condition
\begin{equation}\label{positivity}
f_{r_k}^{<k>} \geq 1\ , \quad k\in \mathbb N
\end{equation}
(which guarantees every $f^{<k>}(z)\in z\mathbb Z_+[[z]])$.
Because $\lim_k r_k=\infty$,
 it follows
that for every $n$ in $\mathbb N$, the sequence
$(f_n^{<k>})_{k=1}^{\infty}$ is eventually constant. Define
$f_n^{<\infty >}= \lim_kf_n^{<k>}$ and $f^{<\infty >}=
\sum_{n=1}^{\infty} f_n^{<\infty >}z^n$. Let $f^{(n)}$ be the
series $f^{<k>}$ where  $k$ is defined by the condition $r_j <n$
if and only if $j<k$.
\newline
Then there is a one-block code $\varphi  : \sigma_{f^{<\infty >}}
\to \sigma_f$ with the following properties.
\begin{enumerate}
\item \label{a} $\varphi $ is  injective. \item \label{aa} For
every $n$, $1-f^{(n+1)}(z) = (1-f^{(n)}(z))/((1-z^n)^{R_n}) \ . $
\item \label{b} $1-f^{<\infty >}(z) =
(1-f(z))/(\prod_{n=1}^{\infty}(1-z^n)^{R_n}) \ . $ \item
\label{bb} If $\varlimsup_n (R_n)^{1/n}< \lambda$, then
$h(\sigma_{f^{<\infty >}}) = \textnormal{log}(\lambda )$. \item
\label{c} If
 $\varlimsup_n (R_n)^{1/n}< \lambda$,
 then $\sigma_f$ is SPR if and only if
$\sigma_{f^{<\infty >}}$ is SPR. \item \label{d} Suppose
\begin{equation}\label{**}
R_n \leq f_n \text{ for all } n, \text{ and also } R_n< f_n \text{
for } n=r_1 \ .
\end{equation}
Then the map $\varphi $ can be chosen to have a magic word.
 \item \label{e}   Suppose
in addition to {\rm (\ref{**})} that $\varlimsup_n (R_n)^{1/n}<
\lambda$. Then the map $\varphi$ can be chosen to have a magic
word and also to be an entropy-conjugacy.
\end{enumerate}
\end{lemma}

\begin{demo}
{\it Construction of $\widehat{\mathcal G}^{<k>}$.} Inductively,
we will define a sequence of labeled graphs
 $\widehat{\mathcal G}^{<k>}$, such that the underlying graph
 $\mathcal G^{<k>}$ is a loop graph for $f^{<k>}$, the
labeling set is the edge set of $\mathcal G^{<1>}$, and the
labeled graph
 $\widehat{\mathcal G}^{<k>}$
defines a one-block code $\varphi _k: \sigma_{f^{<k>}}\to
\sigma_f$ (here, $(\varphi _ky)_i$ is the label of the edge
$y_i$). We let $\mathcal L^{<k>}$ denote the set of (first return)
loops of
 $\mathcal G^{<k>}$, and
we let $\widehat{\mathcal L}^{<k>}$ denote the labeled loops, from
 $\widehat{\mathcal G}^{<k>}$.

If $l$ is a path $l_1\cdots l_j$ of $j$ edges in a labeled graph,
then we denote by $\widehat l$ its label, which is a word
$\widehat l_1\cdots \widehat l_j$, where $\widehat l_i$ is the
label of $l_i$. Similarly when we label a path $l_1\cdots l_j$
with a word $w_1\cdots w_j$, we are defining $\widehat l_i=w_i$.

Choose
  $\mathcal G^{<1>}$
to be a loop graph for $f=f^{<1>}$. Define $\widehat{\mathcal
G}^{<1>}$ by giving each edge a distinct label. Given
  $\widehat{\mathcal G}^{<k>}$,
we construct
  $\widehat{\mathcal G}^{<k+1>}$  from
  $\widehat{\mathcal G}^{<k>}$
as follows.  ($l^0$ will denote the empty word, so $cl^0=c$.)
\begin{itemize}
\item Choose a loop $l^{<k>}=l$ of length $r_k$ in
  $\mathcal G^{<k>}$
(at least one such loop exists by the positivity condition
(\ref{positivity})). \item Now  $\widehat {\mathcal G}^{<k+1>}$ is
the labeled loop graph with exactly the following set of labeled
loops: for every loop $c$ in
  $\mathcal L^{<k>}\setminus \{l\}$, for $n=0, 1,2,3,\dots $,
 put in $\mathcal G^{<k+1>}$ one loop  $\ell (c,n)$,
of length $|c|+nr_k$, with label $\widehat c \widehat l^n$.
\end{itemize}
 The rule
 $\ell (c,n)\mapsto c l^n$ defines a map
 $\mathcal L^{<k+1>}\to
(\mathcal L^{<k>})^*$ which determines a one block code
$\iota_{k+1}: \sigma_{f^{<k+1>}}\to \sigma_{f^{<k>}}$ which
respects the labeling. Therefore the map $\varphi _{k+1}$ is the
composition of $\iota_{k+1},\dots , \iota_2$. Clearly
$f_m^{<k+1>}=f_m^{<k>}-1$ if $m=r_k$, and $f_m^{<k>}=f_m^{<k+1>}$
if $m<r_k$. Each sequence $(f_n^{<k>})_{k=1}^{\infty}$ is
eventually constant and the positivity condition
(\ref{positivity})  does guarantee each $f^{<k>}\in z\mathbb
Z_+[[z]]$. Moreover, the set of labeled loops in
  $\widehat{\mathcal G}^{<k>}$ of length $n$
or less does not change once $r_k>n$. We define
  $\widehat{\mathcal G}^{<\infty >}$ to be the labeled
loop graph whose labeled loops of length $n$ are those in all but
finitely many of the
  $\widehat{\mathcal G}^{<k>}$.
Then
   $\mathcal G^{<\infty >}$
is a loop graph for $f^{<\infty >}$ and the labels of
  $\widehat{\mathcal G}^{<\infty >}$
define a one block code $\varphi : \sigma_{f^{<\infty >}}\to
\sigma_f$.

{\it Injectivity.} For every $k$ in $\mathbb N$, the map
$\iota_{k+1}: \sigma_{f^{<k+1>}}\to \sigma_{f^{<k>}}$ is
injective. (The map $\iota_{k+1}$ is not surjective; a point $w$
of $\sigma_{f^{<k>}}$ fails to be in the image of of $\iota_{k+1}$
if and only if $w$ is forwardly or backwardly asymptotic to the
periodic orbit $\cdots l^{<k>} l^{<k>} l^{<k>}\cdots $ .) It
follows for $1\leq k < \infty$ that
 $\varphi _k: \sigma_{f^{<k>}}\to \sigma_f$ is
injective, with $\text{Image}(\varphi _1)\supset
\text{Image}(\varphi _2)\supset \text{Image}(\varphi _3)\supset
\cdots $ . For $x\in \sigma_f$ and $1\leq k\leq \infty$, define
$V_k(x)$ to be the set of integers $i$ such that there exists $y
\in \sigma_{f^{<k>}}$ such that $\varphi _k(y)=x$ and the edge
$y_i$ begins at the base vertex of $\mathcal G^{<k>}$.
 (Here
$\varphi _{\infty }=\varphi $.) The injectivity of $\varphi _k$
and the structure of $\iota_k$ show for any given $x$ that
$V_1(x)\supset V_2(x)\supset V_3(x) \cdots$ . A point $x$ of
$\sigma_f$ will have a preimage under $\varphi $ in
$\sigma_{f^{<\infty >}}$ if and only if $\cap_{1\leq k<
\infty}V_k(x)$ is unbounded above and below, in which case
\begin{equation}\label{eqV}
V_{\infty }(x)=\bigcap_{1\leq k< \infty} V_k(x)\ .
\end{equation}
Because $\mathcal G^{<\infty>}$ inherits from the $\mathcal
G^{<k>}$'s the property that distinct loops have distinct labels,
this implies that $\varphi $ is injective.

{\it Zeta function and entropy.} For each $n$, if $r_k>n$ then the
number of loops of length $n$ or less in ${\mathcal G}^{<k>}$ is
the same as in $\mathcal{G}^{<\infty>}$. At any finite stage $k$,
$|\mathcal O_n(\sigma_{f^{<k+1>}})| = |\mathcal
O_n(\sigma_{f^{<k>}})|$, except for the $R_n$ values of $k$ at
which $r_k=n$; at each of these values, $|\mathcal
O_n(\sigma_{f^{<k+1>}})| = |\mathcal O_n(\sigma_{f^{<k>}})|-1$.
This description  proves  (\ref{aa})   and also that $|\mathcal
O_n(\sigma_{f^{<\infty >}})| = |\mathcal O_n(\sigma_f)| - R_n$,
for all $n$. Consequently, using
 the product formula for the zeta function,
we also have
\[
1-f^{<\infty >} = \frac{\prod_{n=1}^{\infty}(1-z^n)^{|\mathcal
O_n(\sigma_f)|}} {\prod_{n=1}^{\infty}(1-z^n)^{R_n}} =
\frac{1-f(z)}{\prod_{n=1}^{\infty}(1-z^n)^{R_n}} \
\]
which proves (\ref{b}). By the same product formula,
\[
\Big(\prod_{n=1}^{\infty}(1-z^n)^{R_n}\Big)^{-1} =
\textnormal{exp}\sum_{n=1}^{\infty } \frac 1n Q_n z^n := q(z)
\]
where $ Q_n = \sum_{k|n}kR_k   $
 and therefore
$\varlimsup (Q_n)^{1/n}= \varlimsup (R_n)^{1/n}=
1/\textnormal{rad}(q)$. Finally, given $\varlimsup (R_n)^{1/n}<
\lambda $, the claim (\ref{bb}) follows from
\begin{align*}
h(\sigma_f) &= \varlimsup_n \log |\mathcal O_n(\sigma_f)|^{1/n} =
\varlimsup_n \log\Big(
|\mathcal O_n(\sigma_f)|-R_n\Big)^{1/n} \\
&= \varlimsup_n \log | \mathcal O_n( \sigma_{f^{<\infty
>}})|^{1/n} = h(\sigma_{f^{<\infty >}})\ .
\end{align*}

{\it Strong positive recurrence.} We have $ \text{rad}(f^{<\infty
>}) = \min \{ \text{rad}(f), \text{rad}(q)\} $. If we assume
$\varlimsup_n (R_n)^{1/n}< \lambda $, it follows that
$\text{rad}(q)<1/\lambda$ and
\[
\text{rad}(f^{<\infty >}) > \frac 1{\lambda} \iff \text{rad}(f) >
\frac 1{\lambda} \ .
\]
Therefore $\sigma_{f^{<\infty >}}$ is SPR if and only if
$\sigma_f$ is SPR, proving (\ref{c}).

{\it A magic word.} Consider the inductive construction in terms
of the sequence
   $(\widehat{\mathcal G}^{<k>})$, $k=1,2, \dots $ .
The loops of
  ${\mathcal G}^{<k+1>}$ are loops  $\ell (c,n)$
and the label of  $\ell (c,n)$ in    $\widehat{\mathcal
G}^{<k+1>}$ is $\hat c \hat l^n$.

Assume (\ref{**}) holds. Let $N=r_1$, so, the first deleted loop
$l^{<1>}$ will have length $N$. Pick a loop $W$ of length $N$. At
this initial stage, in
  $\widehat{\mathcal G}^{<1>}$ every loop is labeled
 by its own name, and symbols of $W$ occur in labels
of no loop other than $W$ itself. So, for
  $\widehat{\mathcal G}=
  \widehat{\mathcal G}^{<1>}$, we have the following
\begin{condition} \label{***}
A symbol of $W$ can occur in a label $\widehat l$ of a loop $l$ in
$\widehat{\mathcal G}$ only as part of an initial segment of
$\widehat l$ which equals $W$
\end{condition}
Now suppose that (\ref{***})  holds for $\widehat{\mathcal
G}=\widehat{\mathcal G}^{<k>}$, and the  loop  $l^{<k>}$ has the
property that no symbol of $W$ occurs in $\widehat l^{<k>}$. It
then follows from the construction that
 (\ref{***})
holds for $\widehat{\mathcal G}=\widehat{\mathcal G}^{<k+1>}$.
However, because (\ref{**}) holds, for each $k$ we can  pick for
$l^{<k>}$ a loop which already exists in $\mathcal G^{<1>}$ and
which is not $W$, and which therefore has a label using no symbol
of $W$.

Inductively, then, we construct the $\widehat{\mathcal G}^{<k>}$
to satisfy (\ref{***}), and in the limit obtain $\widehat{\mathcal
G}^{<\infty >}$. The map $\varphi $
 is the one-block code from
$\sigma_{f^{<\infty >}}$ to $\sigma_f$ defined by the labeling of
loops in $\widehat{\mathcal G}^{<\infty >}$. We claim that $W$ is
a magic word for  $\varphi $. To see this, first note that the
labeling of loops in $\widehat{\mathcal G}^{<\infty >}$ inherits
the following properties from the $\widehat{\mathcal G}^{<k>}$:
\begin{enumerate}
\item distinct loops have distinct labels \item a symbol of $W$
 can occur in a label $\widehat l$ of a loop
$l$ in $\widehat{\mathcal G}^{<\infty >}$ only as part of an
initial segment of $\widehat l$ which equals $W$ \item if $WDW$ is
a word in $\sigma_f$, then $WD$ is the label of  a unique
concatenation of loops in $\widehat{\mathcal G}^{<\infty >}$.
\end{enumerate}
Notice that if (3) did not hold, then $\varphi$ would not be
injective, contrary to what we established.

It follows from property (3)
 and the remark
before (\ref{eqV}) that whenever a point $x$ of $\sigma_f$ has the
property that the set $\{n: x[n,n+|W|-1]=W\}$ is unbounded above
and below, then a preimage of $x$ can be constructed in
$\sigma_{f^{<\infty >}}$.

Finally, suppose $x[i,2|W|+|D|-1]=WDW$ and $\varphi  y=x$. It
follows from (2) that both of the edges $y[i]$ and
$y[i+|W|+|D|-1]$
 have as initial vertex the base vertex of the loop graph
$\widehat {\mathcal G}^{<\infty >}$, and therefore
$y[i,|W|+|D|-1]$ is a concatenation of loops in $\widehat
{\mathcal G}^{<\infty >}$. By (3), this concatenation is unique.
Therefore,  $W$ is a magic word for $\varphi$, proving (\ref{d}).

{\it Entropy-conjugacy.} We remark that when $S$ and $T$ are SPR,
it follows already from (\ref{hc-ai-remark}) that the map
$\varphi$ provided by (\ref{d}) is an entropy conjugacy. We now
proceed to the general case, assuming (\ref{**}) and $\varlimsup_n
(R_n)^{1/n}< \lambda $.

Let $\sigma_R$ be the loop shift corresponding to
$R(z)=\sum_{n=1}^{\infty}R_nz^n$. Appealing to (\ref{**}), we
regard $\sigma_R$ as a subsystem of $\sigma_f$. We construct the
map $\varphi$ as above, with the following additional constraint:
for every $k$, the deleted loop $l^{<k>}$ is chosen to be a loop
of $\sigma_R$. As in the argument for (\ref{d}), it follows that
any loop $W$ of $\sigma_f$ which is not a loop of $\sigma_R$ must
be a magic word for $\varphi$.

As $\varphi$ is injective, it is a measure-preserving isomorphism
$(\sigma_{f^{<\infty>}},\nu) \to(\sigma_f,\varphi\nu)$ for any
invariant probability measure $\nu$ of $\sigma_{f^{<\infty>}}$.
Let us consider an ergodic and invariant probability measure $\mu$
for $\sigma_f$ such that  $\varphi^{-1}$ fails to be a.e. defined
(otherwise $\varphi:(\sigma_{f^{<\infty>}},
 \varphi^{-1}\mu)
\to(\sigma_f,\mu)$ would be an isomorphism).
If $W$ is a loop of $\sigma_f$ which is not a loop of $\sigma_R$,
then $\mu\{x:x[0,|W|-1]\}$ must be zero. Therefore $\mu$ is
supported on $\sigma_R$.

It is enough now to prove $h(\sigma_R)<h(\sigma_f)$. Assume first
that $\sigma_R$ is not SPR. Then $h(\sigma_R)= \log \varlimsup_n
R_n^{1/n} < h(\sigma_f)$, proving the claim in this case. Assume
now that $\sigma_R$ is SPR. Pick some $n_0$ such that
$R_{n_0}<f_{n_0}$ and let $R'_n=R_n$ for all $n\ne n_0$ and
$R'_{n_0}=R_{n_0}+1$. Clearly $\varlimsup_n {R'_n}^{1/n}=
\varlimsup_n R_n^{1/n} < h(\sigma_R) \leq h(\sigma_{R'})$. It
follows that $\sigma_{R'}$ is SPR, hence by (\ref{-spr}) of
Proposition \ref{sprprop}, $h(\sigma_R)<h(\sigma_{R'})\leq
h(\sigma_f)$. This implies the claim.

\end{demo}

\section{Main Results}

\begin{lemma}\label{gapprep}
Suppose  $\sigma_F$ and $\sigma_G$ are mixing loop shifts of equal
entropy $\textnormal{log}(\lambda )>0$, and $1\leq \beta < \lambda
$.  Then for any sufficiently large $N\in \mathbb N$, there are
loop shifts $\sigma_f$ and $\sigma_g$ such that
\begin{enumerate}
\item \label{gap0} There are injective one-block codes $\sigma_f
\to \sigma_F$ and $\sigma_g \to \sigma_G$, each of which is an
entropy-conjugacy with a magic word.
 \item
\label{gap1} $|\mathcal O_n(\sigma_f)| = |\mathcal O_n(\sigma_g)|
= 0$ for $n<N$ , \item \label{gap2} $|\mathcal O_n(\sigma_f)| -
|\mathcal O_n(\sigma_g)| = |\mathcal O_n(\sigma_F)| - |\mathcal
O_n(\sigma_G)|\ , \quad n \geq N$, and \item \label{gap3}
$\textnormal{min}\{f_n, g_n\} \geq \beta^n\ , \quad n\geq N\ . $
\end{enumerate}
\end{lemma}

\begin{demo}
Because $\beta < \lambda$  and
$$
\lim_n |\mathcal O_n(\sigma_F)|^{1/n} = \lambda = \lim_n |\mathcal
O_n(\sigma_G)|^{1/n}
$$
for $N$ sufficiently large  we may assume
$$
n\geq N \implies \min_n \{|\mathcal O_n(\sigma_F)|,|\mathcal
O_n(\sigma_G)|\} \geq 2 \lceil \beta^n \rceil \ ,
$$
where $\lceil \beta^n \rceil $ denotes the integer ceiling of
$\beta^n$. Fix such an $N$. Beginning with  $f=F$, we apply Lemma
\ref{lemmazero} repeatedly with $k=z^n$, where $n$ is the length
of a shortest loop, until all loops of length less than $N$ have
been deleted. This produces a series
$$\overline F = \sum_{n=N}^{\infty} \overline F_n z^n
$$
such that  $\sigma_{\overline F}\to \sigma_F$ is an injective
one-block code which  has a magic word and is an entropy-conjugacy
(indeed, $h(\sigma_{z^n})=0$). We also have:
\begin{align*}
1-\overline F &=
(1-F) / \Big( \prod_{1\leq n < N} (1-z^n)^{|\mathcal O_n(\sigma_F)|} \Big) \\
&= (1-F)  \prod_{1\leq n < N}
\Big(\sum_{k=0}^{\infty}z^{nk}\Big)^{|\mathcal O_n(\sigma_F)|} \ .
\end{align*}
Now
 $|\mathcal O_n(\sigma_{\overline F})| = |\mathcal
O_n(\sigma_F)|\geq 2\lceil \beta^n\rceil$ for $n\geq N$, and
$\overline F_n= |\mathcal O_n(  \sigma_{\overline F})|$ for $N\leq
n < 2N$. Let $b_n=\lceil \beta^n\rceil$ and define the series
(polynomial)
$$b(z)=
\sum_{n=N}^{2N-1}b_nz^n \
$$
(so, $\overline F_n -b_n >\beta^n$ for $N\leq n < 2N$). Set
$h=\overline F-b$ and define $f=hb^*=h(1+b+b^2+\cdots )$. By Lemma
\ref{lemmazero}, there is an injective one-block code $\sigma_f\to
\sigma_{\overline F}$ which has a magic word and is
 an entropy-conjugacy.  Indeed,
$h(\sigma_b)<h(\sigma_F)=h(\sigma_f)$ as $b$ defines a finite
state Markov shift whose entropy must increase when adding any of
the remaining loops of $\sigma_F$. Also, the choice of $b$ implies
 $f_n\geq \beta^n$ for $n\geq N$.

Similarly construct $\overline G$ from $G$, and for $b(z)$ exactly
as above set $k =\overline G-b$ and $g=kb^*$. We have proven
(\ref{gap0}), (\ref{gap1}) and (\ref{gap3}). To prove
 (\ref{gap2}), we compute
\begin{align*}
|\mathcal O_n(\sigma_f)| -|\mathcal O_n(\sigma_g)| &=
\big(|\mathcal O_n(\sigma_{\overline F})| -|\mathcal
O_n(\sigma_b)|\big) - \big(|\mathcal O_n(\sigma_{\overline G})|
-|\mathcal
O_n(\sigma_b)|\big) \\
&= |\mathcal O_n(\sigma_{\overline F})| -
|\mathcal O_n(\sigma_{\overline G})| \\
&= |\mathcal O_n(\sigma_F)| -|\mathcal O_n(\sigma_G)| \ ,\quad
n\geq N\ .
\end{align*}
\end{demo}

\begin{theorem}\label{sufficienttheorem}
Suppose $\sigma_F$ and $\sigma_G$ are loop shifts of equal period
and equal entropy $\log (\lambda )$  (possibly infinite), and
\[
\varlimsup_n \big||\mathcal O_n(\sigma_F)| -|\mathcal
O_n(\sigma_G)|\big|^{1/n} < \lambda \ .
\]
Then  $\sigma_F$ and $\sigma_G$ are almost isomorphic, by an
almost isomorphism which induces an entropy-conjugacy of
$\sigma_F$ and $\sigma_G$.
\end{theorem}

\begin{demo}
 Without loss of generality, suppose $\lambda > 1$. To
begin, suppose the period is one, i.e.
  $\sigma_F$ and $\sigma_G$
are mixing. Let $\gamma = \varlimsup \big||\mathcal O_n(\sigma_F)|
-|\mathcal O_n(\sigma_G)|\big|^{1/n}$ . Choose $\beta$ such that
$\gamma < \beta < \lambda $. Choose $f,g$ and $N$  to satisfy the
statement of Lemma \ref{gapprep} for these $\beta$, $\lambda$, $F$
and $G$, with $N$ large enough that
\begin{equation}\label{Nbig}
n\geq N \implies \big||\mathcal O_n(\sigma_F)| -|\mathcal
O_n(\sigma_G)|\big| \leq  \beta^n -1  \ .
\end{equation}
It suffices now to prove
 the theorem with $\sigma_f$ and $\sigma_g$
 in place of
 $\sigma_F$ and $\sigma_G$.

We will apply the Loops Lemma \ref{loopslemma} to $\sigma_f$ and
to $\sigma_g$, using sequences  $(R_n)$ defined for $f$ and for
$g$ as follows:
\begin{align*}
R_n(f) &= \max \{ 0, |\mathcal O_n(\sigma_f)|-|\mathcal O_n(\sigma_g)| \} \ ,\\
R_n(g) &= \max \{ 0, |\mathcal O_n(\sigma_g)|-|\mathcal
O_n(\sigma_f)| \}\ .
\end{align*}
These sequences $(R_n)$ determine the corresponding sequences
$(r_k)$ of the Loops Lemma, as well as Loops Lemma
 sequences of functions
 $f^{(n)}$ and $g^{(n)}$, the corresponding graphs
$\mathcal G^{(n)}(f)$ and $\mathcal G^{(n)}(g)$, and
 the limit functions
$f^{<\infty >}$ and $g^{<\infty >}$. Recall  $f_n^{(n)}$ and
$g_n^{(n)}$ denote the number of (first return) loops of length
$n$ in $\mathcal G^{(n)}(f)$ and $\mathcal G^{(n)}(g)$,
respectively. Using (\ref{aa}) of the Loops Lemma, it is easy to
verify by induction that for every $n$,
\begin{itemize}
\item for $k<n$, the number of loops of length $k$ in $\mathcal
G^{(n)}(f)$ and $\mathcal G^{(n)}(g)$ is the same, and \item
 $f_n^{(n)}-g_n^{(n)} =|\mathcal O_n(\sigma_f)| -
|\mathcal O_n(\sigma_g)|$ .
\end{itemize}

It follows that $f^{<\infty >} = g^{<\infty >}$. Consequently,
there is a map
 $\sigma_{f^{<\infty >}} \to
 \sigma_{g^{<\infty >}}$
which is a one-block code renaming edges, with inverse likewise a
one-block code (here all nonempty words are magic words). By
symmetry, it now remains to show that the  injective  one-block
code $ \sigma_{f^{<\infty >}}\to \sigma_f $ of the Loops Lemma can
be chosen to have a magic word and to be an entropy-conjugacy.

Applying the definition of $R_n(f)$,
 Lemma \ref{gapprep}(\ref{gap2}),
the bound (\ref{Nbig}), and  Lemma  \ref{gapprep}(\ref{gap3}), we
have   for $n\geq N$  that
\begin{align*}
 R_n(f) &=
 \max \{ 0, |\mathcal O_n(\sigma_f)|-|\mathcal O_n(\sigma_g)| \}
= \max \{ 0, |\mathcal O_n(\sigma_F)|-|\mathcal O_n(\sigma_G)| \} \\
&\leq  \beta^n -1 \leq f_n -1\ .
\end{align*}
Moreover, by  Lemma \ref{gapprep}(\ref{gap1}), the length $r_1$ of
the first loop of $\sigma_f$ deleted in the Loops Lemma
construction is equal to $N$. Therefore  conditions (\ref{d}) and
(\ref{e}) of the Loops Lemma are satisfied, and
$\sigma_{f^{<\infty
>}} \to\sigma_f$ can be chosen to be an entropy-conjugacy with a
magic word. By symmetry, we obtain an entropy-conjugacy with a
magic word  $\sigma_f\to\sigma_g$.
This finishes the proof in the mixing case.

Now suppose that
 $\sigma_F$ and $\sigma_G$ have period $p>1$.
Write the $F,G$ in the form
\[
F(z)= \sum_{n=1}^{\infty } a_nz^{pn}\ , \quad G(z)=
\sum_{n=1}^{\infty } b_nz^{pn}\ .
\]
Define $a =\sum a_nz^n$ and $b=\sum b_nz^n$. The loop shifts
$\sigma_a$ and $\sigma_b$ are mixing and by the previous argument
there is a loop shift $\sigma_{c}$, $c=\sum c_nz^n$, with
injective one-block codes $\sigma_c\to \sigma_a$ and $\sigma_c\to
\sigma_b$ with magic words. Let $H(z)=\sum c_nz^{pn}$. It follows
easily that there are injective one-block codes $\sigma_H\to
\sigma_F$ and $\sigma_H\to \sigma_F$ with magic words.
\end{demo}

\begin{theorem}\label{maintheorem}
Suppose $S_A$ and $S_B$ are Markov shifts, and $S_A$ is SPR with
finite entropy. Then the following are equivalent:
\begin{enumerate}
\item \label{3tenors} $S_B$ is SPR, $h(S_A)=h(S_B)$ and
$\textnormal{period}(S_A)= \textnormal{period}(S_B)$.
\item \label{ai} $S_A$ and $S_B$ are almost isomorphic.
 \item\label{hc} $S_B$ is SPR, $S_A$ and $S_B$ are
 entropy-conjugate.
\end{enumerate}
 Moreover, under these conditions there exists an
entropy-conjugacy which is induced by an almost isomorphism.
\end{theorem}

\begin{demo}
(\ref{ai})$\implies$ (\ref{3tenors}). The almost isomorphism
induces a bijection respecting finitary isomorphism between the
ergodic measures of $S_A$ and $S_B$ with full support. Therefore
$S_B$ has a measure of maximal entropy $\mu_B$, and there is a
finitary isomorphism $(S_A,\mu_A)\to (S_B,\mu_B)$. Because
exponential recurrence is an invariant of finitary isomorphism, it
follows from Proposition \ref{sprprop} that $S_B$ is SPR. The
entropy and period of $S_A$ and $S_B$ are already invariants under
 measure-preserving isomorphism of $(S_A,\mu_A)$ and $
(S_B,\mu_B)$.

This last remark also proves (\ref{hc})$\implies$(\ref{3tenors}).

(\ref{3tenors})$\implies$(\ref{ai}) and (\ref{hc}). Pick a vertex
$v$ in the graph with adjacency matrix $A$ which is used to define
the edge shift $S_A$. Because $S_A$ has finite entropy, for each
$n$ the number $F_n$ of first return loops to $v$ is finite. Let
$F=\sum_{n=1}^{\infty }F_nz^n$. The natural injection $\sigma_F\to
S_A$ is a one-block code with a magic word. Similarly choose
$\sigma_G$ for $S_B$. It suffices to prove the theorem with
$\sigma_F$ and $\sigma_G$ in place of $S_A$ and $S_B$.

Let $\text{log}(\lambda )$ denote the entropy of $\sigma_F$ and
$\sigma_G$. Because $S_A$ and $S_B$ are SPR, it follows from
Corollary
 \ref{coro-On} that
\[ \varlimsup \big||\mathcal O_n(\sigma_F)|
-|\mathcal O_n(\sigma_G)|\big|^{1/n} < \lambda \ .
\]
Now Theorem \ref{sufficienttheorem} shows that $\sigma_F$ and
$\sigma_G$ are almost isomorphic, by an almost isomorphism which
induces an entropy-conjugacy.
\end{demo}

Theorem \ref{maintheorem} begs the following questions.
\begin{questions}Is SPR an
 entropy-conjugacy invariant? On the other hand, must positive
recurrent Markov shifts of equal entropy and period be
entropy-conjugate?
\end{questions}

\section{Application to other dynamical systems}\label{application}

Given a measurable system $T$ (as defined in (\ref{measdef})), we
begin by recalling the definition of the  {\it natural extension}
of $T$ (to an invertible system). Let $\overline X$ denote the
subset of $X^{\mathbb Z}$ consisting of the bisequences $(\dots
x_{-1},x_0,x_1,\dots)$ such that $T(x_i)=x_{i+1}$ for all $i$.
Associated to $\overline X$ is the $\sigma$-algebra generated by
requiring every coordinate projection to be measurable. Let
$\overline T$ be the bijection $\overline X \to \overline X$
defined by applying $T$ coordinatewise. The projection of
$\overline X$ onto the first coordinate defines a map which
intertwines $T$ and $\overline T$ (if $T$ is invertible, then this
coordinate projection is a bijection). This coordinate projection
induces an entropy-respecting bijection $\mathcal M(\overline
T)\to \mathcal M(T)$ of invariant Borel probabilities. Often
measurable phenomena are more conveniently studied on invertible
systems, and the natural extension is used for this.

Now supose $\gamma$ is an entropy conjugacy of invertible systems.
We say  $\gamma$ {\it refines the hyperbolic structure} if the
entropy-negligible set of Definition \ref{entcondef} can be chosen
such that for $x,y$ in its complement,
\begin{align*}
 \lim_{n\to\infty} \text{dist}(\gamma (T^nx),\gamma(T^ny))=0
 &\implies \lim_{n\to\infty} \text{dist}(T^nx,T^ny)=0 \ , \quad
    \text{and} \\
 \lim_{n\to-\infty} \text{dist}(\gamma (T^nx),\gamma(T^ny))=0
 &\implies \lim_{n\to-\infty} \text{dist}(T^nx,T^ny)=0 \ .
\end{align*}
We say $\gamma$ {\it preserves the hyperbolic structure} if both
$\gamma$ and its inverse refine it.

Our next theorem is essentially an assembly of results from
 \cite{SIM,AFF,EE,QFT}, collected for the statement of
Theorem \ref{piecewiseclassif}. To prepare,  recall some
definitions. A {\it piecewise monotonic} interval map  is
 a map $f:[0,1]\to[0,1]$ such that there is a finite partition of $[0,1]$ into
intervals on each of which the restriction of $f$ is continuous
and strictly monotonic. The map $f$ is said to be {\it completely
discontinuous} if the previously mentioned partition can be chosen
so that all its endpoints except $0$ and $1$ are discontinuities
of $f$.

A {\it multi-dimensional $\beta$-transformation} \cite{AFF} is a
map defined by $x\mapsto B.x \mod \mathbb Z^d$ for $x\in[0,1]^d$,
$d=1,2,\dots$ and $B$ an expanding affine map of $\mathbb R^d$.

We will not define the class of {\it subshifts of quasi finite
type}  \cite{QFT};
 this is  a useful class of
finite state subshifts which contains the class of subshifts of
finite type, but is much larger (for example, every finite
positive entropy is achieved by a subshift of quasi finite type).

\begin{theorem}\label{classlist} The following measurable dynamical systems
have natural extensions which are entropy-conjugate to the
disjoint union of finitely many SPR Markov shifts of equal
entropy.
\begin{enumerate}
 \item Subshifts of quasi finite type \cite{QFT}.
 \item Piecewise monotonic interval maps with non-zero topological
entropy.
 \item The multi-dimensional $\beta$-transformations \cite{AFF}.
 \item  $C^\infty$ smooth entropy-expanding maps including, e.g., smooth interval maps
 with non-zero topological entropy \cite{SIM} or $(x,y)\mapsto(a-x^2+\eps y,b-y^2+\eps x)$
 with $a,b<2$ close to $2$ and $|\eps|$ sufficiently small \cite{BSMF,EE}.
 \end{enumerate}
Moreover, the entropy-conjugacy from the natural extension to the
SPR Markov shift can be chosen to refine the hyperbolic structure.
In cases (1) and (3) and for completely discontinuous piecewise
monotonic maps, the entropy-conjugacy can be chosen to preserve
the hyperbolic structure.\footnote{In cases (2) and (4), it is the
hyperbolic structure of the symbolic dynamics which is preserved.}
\end{theorem}
\begin{proof}[Proof sketch]
(1) follows from \cite[Theorem 3]{QFT}.

(2) follows from \cite[Lemma 2]{QFT}.

(3) also follows from \cite[Lemma 2]{QFT}.

(4)
 An  argument similar  to that in \cite[Lemma 2]{QFT} can be
applied to $C^\infty$ smooth entropy-expanding maps in arbitrary
dimension (see \cite{EE}). More specifically, in the proof in
section III of \cite{EE}, the condition that $\hat\mu(\Sigma_*)$
is zero can be replaced by the condition that it is very small.
Then we can conclude  that the Markov shift is SPR as in
\cite{QFT}, Corollary at the end of section ``Entropy at
infinity''.

 The proof of the preservation of the hyperbolic
structure is done in two steps.  To begin, we recall that
multidimensional $\beta$-transformations and completely
discontinuous piecewise monotonic maps have natural extensions
that are topologically conjugate to those of their symbolic
dynamics. These are subshifts of quasi-finite type (or Q.F.T.) as
shown in \cite{QFT}. So we can restrict to  Q.F.T.

We use terminology and facts from \cite{QFT}. A Markov shift (the
so-called complete Markov diagram) is built over the Q.F.T. Let
$\gamma$ be the entropy-conjugacy  of the original subshift into
that Markov shift. $\gamma$ obviously refines the hyperbolic
structure as the canonical partition of the Markov shift is finer
than the original canonical partition, using the identification
$\gamma$.

On the other hand, for $x$ in the Markov shift, $x_0$ is
determined by the sequence $(\gamma^{-1} x)_n, n\leq0$. Thus the
unstable relation is preserved by $\gamma$.

One cannot replace $n\leq0$ by $n\geq0$ here. But we have the
following property for a.e. point $x$ in the Markov shift: $x_n$
is determined by the finite word $(\gamma x)_k, {n-\ell}\leq k\leq
n$ with $\ell=\ell(\sigma^n x)$. By stationarity
$\ell(\sigma^nx)\not\to\infty$ for a.e. $x$. Hence the future in
the subshift determines $x_{n_0}$ for some $n_0$. But $x_{n_0}$
and the sequence $(\gamma^{-1} x)_k, k\geq n_0$ determine
$x_k,k\geq n_0$. Thus the stable relation is also preserved by
$\gamma$. This concludes the proof for the case of Q.F.T.
\end{proof}

Combining  Theorem \ref{classlist} with our main result, Theorem
\ref{maintheorem}, we obtain

\begin{theorem}\label{piecewiseclassif}
A system in any of the classes listed in Theorem   \ref{classlist}
has finite entropy and only finitely many ergodic measures of
maximal entropy. Two systems from the list have natural extensions
which are entropy conjugate if and only if  they have equal
entropy and for each $p\in \mathbb N$ the same number of
 ergodic maximal entropy measures of period $p$.
\end{theorem}

\begin{corollary}
Two topologically  mixing piecewise monotonic interval maps are
entropy conjugate if and only if they have the same entropy.
\end{corollary}

\begin{proof}
A topologically mixing piecewise monotonic map has a unique
measure of maximal entropy, and this measure is mixing \cite{H2}.
 So, in this case the invariant above simplifies to the
entropy of that measure, which is the topological entropy by the
variational principle \cite{Walters-book}.
\end{proof}

\section{Remarks} \label{remarkssection}

We continue our convention that ``Markov shift'' means
``irreducible Markov shift'', even in the finite state case.

\begin{remark}[Magic words]
Magic words come from the finite state coding theory of Markov
shifts \cite{Ki,LM}. They play a basic role in the structure
theory of finite-to-one codes between finite-state Markov shifts
(where they are
 defined more generally for maps which are finite-to-one
but not necessarily one-to-one almost everywhere). In the past, a
major step in some code constructions was to show choices could be
made to guarantee existence of a magic word, so that the resulting
code would be one-to-one almost everywhere. (For example, the
construction of a magic word
 is the essential difference between \cite{AM} and \cite{Parry}.)
Ashley's Replacement Theorem \cite{Ki,LM} has largely eliminated
this step for codes of finite state systems.
\end{remark}

\begin{remark}[Almost conjugacy] \label{adlermarcus}
 \cite{AM,Ki,LM}
Two finite state Markov shifts $S,T$ are {\it almost conjugate}
\cite{LM}, equivalently {\it almost topologically conjugate}
\cite{AM,Ki}, if there is another Markov shift $R$ and surjective
block codes $R\to S$, $R\to T$ each of which are one-to-one almost
everywhere. The Adler-Marcus Theorem \cite{AM} is that entropy and
period are complete invariants of almost topological conjugacy for
finite state Markov shifts. (Here each of the
 maps $R\to S$, $R\to T$
must have a magic word, and the content  of ``a.e.'' contains no
more information
 than  we have  in Propositions \ref{expspeed} and \ref{magicprop}.)

The classification of countable state Markov shifts up to almost
conjugacy is extremely difficult and remains an open problem
\cite{DF}. Doris Fiebig has the most general result: for Markov
shifts $S,T$ of equal finite positive entropy, there is a Markov
shift $R$ of equal entropy and countable-to-one biclosing
continuous surjections $R\to S, R\to T$. Here ``biclosing'' is a
bonus and it is impossible to replace ``countable-to-one'' with
 ``uniformly finite-to-one''   \cite{DF}.
Note, this result places no
restriction on the relative classes of $S$ and $T$ (which can be
transient, null recurrent, positive recurrent or strong positive
recurrent, each independent of the other).
\end{remark}

\begin{remark}[Almost isomorphism]
Our definition of ``almost isomorphism'' is analogous to the
definition of ``almost conjugacy''; for countable state Markov
shifts, it seems more tractable and  perhaps better suited to
studying invariant measures. Injective block codes between finite
state Markov shifts of equal entropy must be homeomorphisms, so
``almost isomorphism'' does not  appear as a meaningful distinct
relation in the finite state category.
\end{remark}

\begin{remark}[Other Markov measures]\label{othermarkov}
The Adler-Marcus Theorem gave good finitary codes between finite
state Markov shifts of equal entropy and period. By a Markov
chain, we will mean a Markov shift together with an invariant
Markov measure. The subsequent search for good finitary codes
between  finite state Markov chains has a sizeable history (see
\cite{MT,MoTu,Pa1,PS,S1,S2,S3} and their references); for years,
there was an attempt to produce general codes by way of an almost
topological conjugacy.
 Marcus and Tuncel destroyed these hopes with their weight-per-symbol
polytope invariants \cite{MT}. The paper \cite{MoTu} of Mouat and
Tuncel effectively exploits
 one alternate approach, using
ordered algebra. The ``loops method'' of this paper is another
approach; we hope it will be of  use for Markov chains beyond the
case of the measure of maximal entropy.
 Both the Mouat-Tuncel construction and the
Loops Lemma method yield magic word isomorphisms when they work.
We expect that the known invariants (ratio group \cite{Kri},
weights group with distinguished coset \cite{PS}, beta function
\cite{S1}) of hyperbolic-structure-preserving isomorphism
\cite{S2,S3} with finite-expected coding time will be the only
obstructions to the existence of a magic word isomorphism between
finite state mixing Markov chains.
\end{remark}

\begin{remark}[Local compactness and loop shifts]
Where we use ``irreducible Markov shift'', Gurevich and Savchenko
\cite{GS} use ``indecomposable symbolic Markov chain''. Their use
of ``symbolic'' rather than ``topological'' was to emphasize that
topological issues are not always paramount.  The following simple
construction illustrates this point in our context.

Let $\sigma_f$ be any loop shift such that $f(1)=\infty$ (i.e.,
$f$ is not just a polynomial). Then the space $\sigma_f$ will be
neither locally compact nor $\sigma$-compact. The loop shift is
still closely related to a locally compact system. The graph in
Figure 1 describes the construction of a locally compact $S_A$
from $f$, such that the natural injective one block code
$\sigma_f\to S_A$ has many magic words and induces a bijection of
spaces of invariant Borel probabilities, $\mathcal M(\sigma_f)\to
\mathcal M(S_A)$.

\bigbreak

\begin{center}
\epsfig{figure=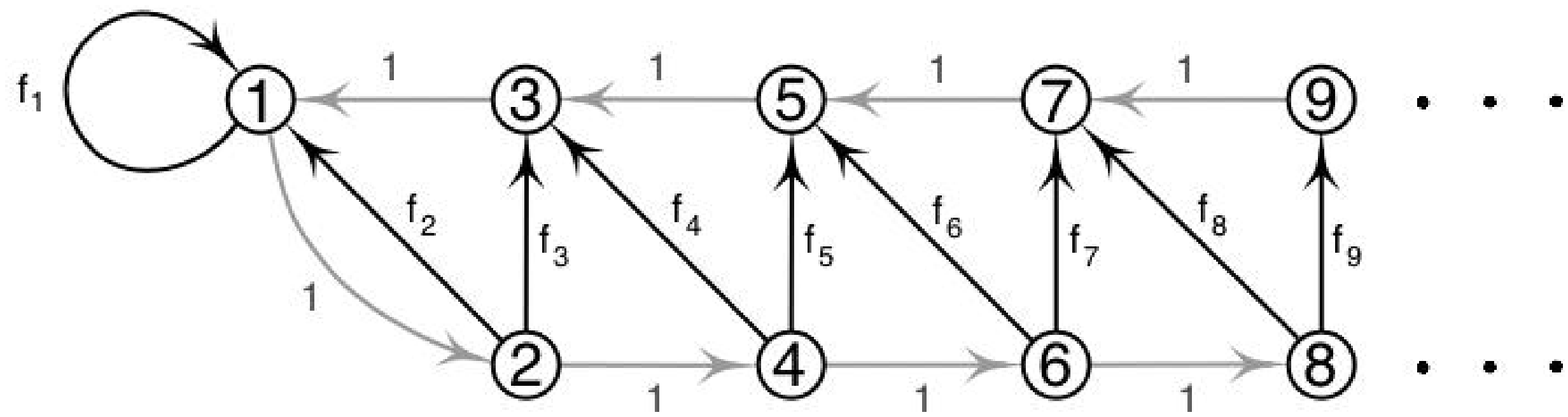,height=3cm}

{\small {\sl Figure 1.}}
\end{center}

\end{remark}


\begin{thebibliography}{99999}

\bibitem{AM} R. Adler and B. Marcus,
    {\em Topological entropy and equivalence of dynamical
    systems}, Memoirs Amer. Math. Soc. {\bf 20} (1979), no. 219.


\bibitem{BGMY} L. Block, J. Guckenheimer, M. Misiurewicz and L.S. Young,
    {\em Periodic points and topological entropy of one-dimensional maps},
    Global theory of dynamical systems, pp. 18--34,
    Lecture Notes in Math. {\bf 819} (1980), Springer-Verlag.

\bibitem{BowenLNM} R. Bowen,
    {\em Equilibrium states and the ergodic theory of Anosov diffeomorphisms},
    Lecture Notes in Math. {\bf 470} (1975), Springer-Verlag.




\bibitem{BBG} M. Boyle, J. Buzzi and R. Gomez,
    {\em Good potentials for almost isomorphism of countable state
    Markov shifts}, preliminary manuscript, 2004.

\bibitem{SIM} J. Buzzi,
    {\em Intrinsic ergodicity of smooth interval maps},
    Israel J. Math. {\bf 100} (1997), 125--161.

\bibitem{AFF} J. Buzzi,
    {\em Intrinsic ergodicity of affine maps on $[0,1]^d$},
    Monat. Math. {\bf 124} (1997), 97-118.

\bibitem{BSMF} J. Buzzi,
    {\em Ergodicit\'e intrins\`eques de produits fibr\'es d'applications
    chaotiques unidimensionnelles},
    Bull. Soc. Math. France {\bf 126} (1998), 51-77; English
    preprint version at {\tt
    http://math.polytechnique.fr/cmat/buzzi}.

\bibitem{EE} J. Buzzi, On entropy-expanding maps,
    {\em preprint CMAT}, 2000.

\bibitem{QFT}
J. Buzzi, Subshifts of quasi-finite type,
    {\em preprint CMAT,} 2003.

\bibitem{BS}  J. Buzzi, O. Sarig,
    {\em Uniqueness Of Equilibrium Measures For Countable
    Markov Shifts And Multidimensional Piecewise Expanding Maps},
    Ergod. Th. \& Dynam. Syst. {\bf 23} (2003), 1383--1400.

\bibitem{EEC} E.E.Cummings,
    {\em Collected works 1913-1962},
    Harcourt Brace Jovanovich, Inc. (1972).

\bibitem{DF} D. Fiebig,
    {\em Common extensions for locally compact Markov shifts},
    Monatsh. Math. {\bf 132} (2001), no. 4, 289--301.

\bibitem{DF2} D. Fiebig,
    {\em Graphs with pre-assigned Salama entropies and optimal
    degrees  for locally compact Markov shifts},
     Ergodic Theory Dynam. Systems {\bf 23} (2003), no. 4,
    1093--1124.

\bibitem{DF3} D. Fiebig,
    {\em Factor theorems for locally compact Markov shifts}
    Forum Math. {\bf 14} (2002), no. 4, 623--640.

\bibitem{DF4} D. Fiebig,
    {\em Factor maps, entropy and fiber cardinality for Markov
    shifts},
     Rocky Mountain J. Math. {\bf 31} (2001), no. 3, 955--986.

\bibitem{FF} D. Fiebig and U. Fiebig,
    {\em Entropy and finite generators for locally compact
    subshifts},
     Ergodic Theory Dynam. Systems {\bf 17} (1997), no. 2,
    349--368.

\bibitem{FFY} D. Fiebig, U.-F. Fiebig and M. Yuri,
    {\em Pressure and equilibrium states for countable
    state Markov shifts},
    Israel J. Math. 131 (2002), 221--257.

\bibitem{UF2} U.-R. Fiebig,
    {\em A return time invariant for finitary isomorphisms},
    Ergodic Theory Dynam. Systems {\bf 4} (1984), no. 2, 225--231.

\bibitem{UF} U.-R. Fiebig,
    {\em Symbolic dynamics and locally compact Markov shifts},
    1996. Habilitationsschrift, U. Heidelberg.


\bibitem{G} R. Gomez,
        {\em Positive K-theory for finitary isomorphisms of Markov chains},
    Ergodic Theory  Dynam. Systems {\bf 23} (2003), pp 1485-1504.

\bibitem{G0} B. M. Gurevich,
    {\em Shift entropy and Markov measures in the path space of
    a denumerable graph} (Russian),  Dokl. Akad. Nauk SSSR
    {\bf 187} (1969), 715--718; English translation: Soviet Math.
    Dolk. {\bf 10}, 4, 911--915.

\bibitem{G1} B. M. Gurevich,
    {\em Stably recurrent nonnegative matrices} (Russian),
    Uspekhi Mat. Nauk {\bf  51}  (1996),  no. 3(309), 195--196;
    translation in  Russian Math. Surveys {\bf 51}  (1996),
    no. 3, 551--552.


\bibitem{GP} B. M. Gurevich and A.B. Polyakov,
    {\em Stably recurrent functions on the loop space of a countable graph},
    Uspekhi Mat. Nauk {\bf 54} (1999), no. 6(330), 157--158;
    translation in Russian Math. Surveys {\bf 54} (1999), no. 6, 1242--1243.

\bibitem{GS} B. M. Gurevich and S. Savchenko,
    {\em Thermodynamical formalism for symbolic Markov
    chains with a countable number of states} (Russian),
    Uspekhi Mat. Nauk {\bf 53} (1998),  3--106;
    translation in  Russian Math. Surveys {\bf 53} (1998), 245--344.

\bibitem{H1} F. Hofbauer,
    {\em $\beta $-shifts have unique maximal measure},
     Monatsh. Math. {\bf 85}  (1978), no. 3, 189--198.

\bibitem{H2} F. Hofbauer,
    {\em On intrinsic ergodicity of piecewise monotonic
    transformations with  positive entropy I},
     Israel J. Math. {\bf 34} (1979), no. 3, 213--237;
    II Israel J. Math. {\bf 38}  (1981), no. 1-2, 107--115.

\bibitem{H3} F. Hofbauer,
    {\em Piecewise invertible dynamical systems},
    Probab. Theory Relat. Fields  {\bf 72}  (1986),  no. 3, 359--386.

\bibitem{KS} M. Keane and M. Smorodinsky,
        {\em Finitary isomorphisms of irreducible Markov shifts},
         Israel J. Math. {\bf 34} (1979), no. 4, 281--286 (1980).

\bibitem{K1} G. Keller,
    {\em Markov extensions, zeta functions, and Fredholm theory
    for piecewise invertible dynamical systems},
    Trans. Amer. Math. Soc.  {\bf 314}  (1989),  no. 2, 433--497.

\bibitem{K2} G. Keller,
    {\em Lifting measures to Markov extensions},
    Monatsh. Math.  {\bf 108}  (1989),  no. 2-3, 183--200.

\bibitem{Ki} B.P. Kitchens,
        {\em Symbolic dynamics. One-sided, two-sided
        and countable state Markov shifts},
        Springer-Verlag (1998).

\bibitem{Kre} U. Krengel,
 {\em Ergodic theorems,}
 Walter de Gruyter, (1985).

\bibitem{Kri} W. Krieger,
    {\em On the finitary isomorphisms of Markov shifts that have
    finite expected coding time},
    Z. Wahrsch. Verw. Gebiete {\bf 65} (1983), no. 2, 323--328.

\bibitem{Ledrappier} F. Ledrappier,
    {\em Principe variationnel et syst\`emes dynamiques
    symboliques},
    Z. Wahrsch. Verw. Gebiete {\bf 30} (1974),  185--202.

\bibitem{LM} D. Lind and  B. Marcus,
        {\em An introduction to Symbolic Dynamics and Coding},
        Cambridge University Press (1995).

\bibitem{LSV} C. Liverani, B. Saussol, S. Vaienti,
    {\em Conformal measure and decay of correlation for covering weighted
    systems,}
   Ergodic Theory Dynam. Systems {\bf 18} (1998), 1399--1420.

\bibitem{MT} B. Marcus and S. Tuncel,
    {\em The weight-per-symbol polytope and scaffolds of
    invariants associated with Markov chains},
     Ergodic Theory Dynam. Systems {\bf 11} (1991),
    129--180.

\bibitem{MU1} R. D. Mauldin and M. Urba\'nski,
    Gibbs states on the symbolic space over an infinite alphabet.
    Israel J. Math. {\bf 125} (2001), 93--130.

\bibitem{MoTu} R. Mouat and S. Tuncel,
    {\em Constructing finitary isomorphisms with finite
    expected coding times},
    Israel J. Math. {\bf 132} (2002), 359--372.

\bibitem{Orn} D. Ornstein,
    {\em Bernoulli shifts with the same entropy are isomorphic.}
    Advances in Math. {\bf 4} (1970), 337--352 (1970).

\bibitem{Parry} W. Parry,
    {\em A finitary classification of topological Markov chains
    and sofic systems},
    Bull. London Math. Soc. {\bf 9} (1977), no. 1, 86-92.

\bibitem{Pa1} W. Parry,
    {\em Notes on coding problems for finite state processes},
    Bull. London Math. Soc.{\bf 23} (1991), 1-33.


\bibitem{PS} W. Parry and K. Schmidt,
    {\em Natural coefficients and invariants for Markov shifts},
    Inventiones Math. {\bf 76} (1984), 15-32.


\bibitem{RudMarkov} D. Rudolph,
    {\em A mixing Markov chain with exponentially decaying return
    times is finitarily Bernoulli},
     Ergodic Theory Dynam. Systems {\bf 2}
    (1982), no. 1, 85--97.

\bibitem{ruette} S. Ruette,
    {\em On the Vere-Jones classification and existence of
    maximal measures for countable topological Markov chains},
    Pacific J. Math. {\bf 209} (2003), 365--380.

\bibitem{salama1} I. Salama,
    {\em Topological entropy and recurrence of countable chains},
     Pacific J. Math. {\bf  134}  (1988),  no. 2, 325--341;
    corrections {\bf 140}  (1989),  no. 2, 397.

\bibitem{salama2} I. Salama,
    {\em On the recurrence of countable topological Markov chains},
    in Symbolic dynamics and its applications (New Haven, CT, 1991),
     Contemp. Math., {\bf 135} (1992), 349--360,
    Amer. Math. Soc., Providence, RI,

\bibitem{Sarig0} O. Sarig,
    {\em Thermodynamic formalism for countable Markov shifts},
    Ergodic Theory Dynam. Systems {\bf 19} (1999), no. 6, 1565--1593.

\bibitem{Sarig}
    O. Sarig,
    {\em Thermodynamic formalism for null recurrent potentials},
    Israel J. Math. {\bf 121} (2001), 285--311.

\bibitem{Sarig2}
    O. Sarig,
    {\em Phase transitions for countable Markov shifts},
    Comm. Math. Phys. {\bf  217}  (2001),  no. 3, 555--577.

\bibitem{Sarig4}
    O. Sarig,
    {\em Subexponential decay of correlations,}
     Invent. Math. {\bf 150} (2002), 629-653.

\bibitem{Sarig3}
    O. Sarig,
    {\em Existence of Gibbs measures for countable Markov shifts},
    Proc. Amer. Math. Soc. {\bf  131}  (2003),  no. 6, 1751--1758

\bibitem{S1} K. Schmidt,
    {\em Invariants for finitary isomorphisms with finite
    expected code lengths},
    Invent. Math. {\bf 76} (1984), no. 1, 33--40.

\bibitem{S2} K. Schmidt,
    {\em Hyperbolic structure preserving isomorphisms of Markov
    shifts},
    Israel J. Math. {\bf 58} (1987), no. 2, 225--242.

\bibitem{S3} K. Schmidt.
    {\em Hyperbolic structure preserving isomorphisms of Markov
    shifts},
    Israel J. Math. {\bf 55} (1986), no. 2, 213--228.

\bibitem{Tak} Y. Takahashi,
    {\em Isomorphisms of $\beta $-automorphisms to Markov
    automorphisms},
    Osaka J. Math. {\bf 10} (1973), 175--184.

\bibitem{Tak2} Y. Takahashi,
    {\em Shift with orbit basis and realization of
    one dimensional maps},
    Osaka J. Math. {\bf 20} (1983), 599-629;
    Correction  Osaka J. Math. {\bf 22} (1985), 637.

\bibitem{Thaler} M. Thaler,
   {\em Estimates of the invariant
   densities of endomorphisms with indifferent fixed points,}
    Israel J. Math. {\bf 37}
   (1980), 303--314.

\bibitem{MU2} M. Urba\'nski and A.  Zdunik,
    {\em Hausdorff dimension of harmonic measure for self-conformal
    sets},
    Adv. Math. {\bf 171} (2002), no. 1, 1--58.


\bibitem{Vere-Jones1}
    D. Vere-Jones,
    {\em Ergodic properties of non-negative matrices},
    Pacific J. Math. {\bf 22} (1967), 361--386.

\bibitem{Vere-Jones2}
    D. Vere-Jones,
    {\em Geometric ergodicity in denumerable Markov chains},
    Quart. J. Math. Oxford (2) {\bf 13} (1962), 7--28.

\bibitem{Walters0}
    P. Walters,
    {\em Ruelle's operator theorem and $g$-measures},
    Trans. Amer. Math. Soc. {\bf 214} (1975), 375--387.


\bibitem{Walters-book}
    P. Walters,
    {\it An introduction to ergodic theory},
    Springer-Verlag (1982).



\bibitem{Young} L.-S. Young,
    {\em Recurrence times and rates of mixing},
    Israel J. Math. {\bf 110} (1999), 153--188.

\bibitem{Yuri0} M. Yuri,
    {\em On the convergence to equilibrium states of certain
    non-hyperbolic systems},
    Ergod. Th. \& Dynam. Syst. {\bf 17} (1997), 977-1000.

\end{thebibliography}
\end{document}